\documentclass{elsart}%
\usepackage{amsfonts}
\usepackage{amsmath}
\usepackage{amssymb}
\usepackage{graphicx}%
\journal{Topology}
\begin{document}
%
\begin{frontmatter}%
%

\title{A maximum principle for combinatorial Yamabe flow}%
%

\author{David Glickenstein}%
%

\address
{Department of Mathematics, University of Arizona, Tucson, AZ 85721, USA}%
%

\begin{abstract}%
This article studies a discrete geometric structure on triangulated manifolds
and an associated curvature flow (combinatorial Yamabe flow). The associated
evolution of curvature appears to be like a heat equation on graphs, but it
can be shown to not satisfy the maximum principle. The notion of a
parabolic-like operator is introduced as an operator which satisfies the
maximum principle, but may not be parabolic in the usual sense of operators on
graphs. A maximum principle is derived for the curvature of combinatorial
Yamabe flow under certain assumptions on the triangulation, and hence the heat
operator is shown to be parabolic-like. The maximum principle then allows a
characterization of the curvature as well was a proof of long term existence
of the flow.
\end{abstract}%
%

%

\end{frontmatter}%

\section{Introduction}

In \cite{glickensteincombinatorialyamabeflow} we introduced the combinatorial
Yamabe flow on three-dimensional piecewise linear complexes as an analogue of
the smooth Yamabe flow (see \cite{hamilton:surfaces}, \cite{ye} for the Yamabe
flow and \cite{leeparker} for a look at the Yamabe problem). The complexes are
geometric in the sense that each Euclidean tetrahedron is given a metric
structure by having the edge lengths determined by weights $r_{i}$ defined at
each vertex $i$. The length of an edge $\left\{  i,j\right\}  $ is defined to
be $r_{i}+r_{j};$ all such structures are called a conformal class as they are
a conformal deformation of the triangulation where all edges are the same
length. For each (nondegenerate) tetrahedron, we can think of the structure as
coming from four mutually tangent spheres such that the centers are connected
by the edges of the tetrahedron and the $i$th sphere has radius $r_{i}$. We
showed that under sufficient long term existence conditions, the flow
converges to constant curvature.

The three-dimensional combinatorial Yamabe flow was inspired by the work of
Chow and Luo \cite{chowluo} (see also \cite{luoyamabe}). They looked at a
combinatorial Ricci flow on two dimensional simplicial complexes and showed
that the equation satisfies a maximum principle. The maximum principle, which
says that the maximum of a solution decreases and the minimum increases, is
one of the most useful concepts in the study of the heat equation and other
parabolic partial differential equations. Maximum principle techniques have
been used to great benefit in the smooth category. We are especially inspired
by Hamilton's work on the Ricci flow (see \cite{hamilton:singularities} and
\cite{caochow}). It is in general very difficult to prove a maximum principle,
or even to prove short term existence of solutions, when you do not have a
strictly parabolic equation.

In this paper we investigate an analytic result about the flow which leads to
a long term existence result for some structures. We find that the evolution
of curvature admits a maximum principle under certain assumptions on the
triangulation. It is especially interesting that we are able to derive a
maximum principle even in a situation where the evolution is not parabolic in
the usual sense of graph Laplacians. We thus introduce the notion of
parabolic-like operators which satisfy the maximum principle for a given
function. We can then show that under sufficient assumptions the evolution of
curvature is parabolic-like.

\section{Parabolic-like operators}

The weighted (unnormalized) Laplacian on a graph $G=\left(  V,E\right)  ,$
where $V$ are the vertices and $E$ are the edges, is defined as the operator
\[
\left(  \triangle f\right)  _{i}=\sum_{\left\{  i,j\right\}  \in E}%
a_{ij}\left(  f_{j}-f_{i}\right)
\]
for each $i\in V,$ where the coefficients satisfy $a_{ij}=a_{ji}$ and
$a_{ij}\geq0$ (see, for instance, \cite{chung}). The operator $\triangle$
takes functions on $V$ to functions on $V.$ The coefficients depend on the
edge $\left\{  i,j\right\}  $ (compare \cite{chungvertex}). The symmetry
condition is simply self-adjointness with respect to the Euclidean metric, so
we can replace it with another self-adjointness condition. That is, we can
define an inner product $\left\langle f,g\right\rangle _{b}=\sum_{i}b_{i}%
f_{i}g_{i}$ if we are given coefficients $\left\{  b_{i}\right\}  .$ An
operator $S$ is self-adjoint with respect to $b$ if%
\[
\left\langle Sf,g\right\rangle _{b}=\left\langle f,Sg\right\rangle _{b}.
\]
It is clear that symmetry corresponds to being self-adjoint with respect to
the inner product determined by $b_{i}=1$ for all $i.$ We shall call $\left\{
b_{i}\right\}  $ a (positive definite) metric if $b_{i}>0$ for all $i.$ In
order to match with the notation in \cite{glickensteincombinatorialyamabeflow}%
, we shall let $\mathcal{S}_{0}$ denote the vertices and $\mathcal{S}_{1}$
denote the edges. In later sections we will usually consider the inner product
coming from the metric $\left\{  r_{i}\right\}  _{i\in\mathcal{S}_{0}}.$

We define a (discrete) parabolic operator on functions
\begin{equation}
f:\mathcal{S}_{0}\times\left[  A,Z\right)  \rightarrow\mathbb{R}%
,\label{functiontypes}%
\end{equation}
where $\mathcal{S}_{0}$ is a discrete set of vertices and $\left[  A,Z\right)
\subset\mathbb{R}$, as follows. First we shall call $\mathcal{F}$ the class of
functions of the form (\ref{functiontypes}). We write the evaluation of $f$ at
the point $\left(  i,t\right)  $ as $f_{i}\left(  t\right)  .$

\begin{defn}
\label{parabolic def}An operator
\[
P:\mathcal{F}\rightarrow\mathcal{F}%
\]
of the form%
\[
\left(  Pf\right)  _{i}\left(  t\right)  =\frac{df_{i}\left(  t\right)  }%
{dt}-\sum_{\left\{  i,j\right\}  \in\mathcal{S}_{1}}a_{ij}\left(  t\right)
\left(  f_{j}\left(  t\right)  -f_{i}\left(  t\right)  \right)  ,
\]
where $a_{ij}:\left[  A,Z\right)  \rightarrow\mathbb{R}$ are self-adjoint with
respect to some metrics $\left\{  b_{i}\left(  t\right)  \right\}
_{i\in\mathcal{S}_{0}},$ is called \emph{parabolic }if $a_{ij}\left(
t\right)  \geq0$ for all $i,j\in\mathcal{S}_{0}$ and for all $t\in\left[
A,Z\right)  .$
\end{defn}

Parabolic operators are of the form $\frac{d}{dt}-\triangle$ if $\triangle$ is
defined to be an appropriate Laplacian with weights. We note that it is easy
to prove a maximum principle for these operators as follows:

\begin{prop}
If $P$ is parabolic and $f$ is a solution to $Pf=0,$ then $f$ satisfies%
\begin{align*}
\frac{df_{M}}{dt}\left(  t\right)   &  \leq0\\
\frac{df_{m}}{dt}\left(  t\right)   &  \geq0
\end{align*}
where $M,m\in\mathcal{S}_{0}$ such that%
\begin{align*}
f_{M}\left(  t\right)   &  \doteqdot\max_{i\in\mathcal{S}_{0}}\left\{
f_{i}\left(  t\right)  \right\} \\
f_{m}\left(  t\right)   &  \doteqdot\min_{i\in\mathcal{S}_{0}}\left\{
f_{i}\left(  t\right)  \right\}  .
\end{align*}

\end{prop}

\begin{pf}
Since $Pf=0$ we have, for a given $t,$
\[
\frac{df_{M}}{dt}=\sum_{\left\{  M,j\right\}  \in\mathcal{S}_{1}}a_{Mj}\left(
t\right)  \left(  f_{j}\left(  t\right)  -f_{M}\left(  t\right)  \right)
\]
and since $a_{Mj}\geq0$ and $f_{j}\left(  t\right)  \leq f_{M}\left(
t\right)  $ for all $i\in\mathcal{S}_{0}-\left\{  M\right\}  $ we see that
\[
\frac{df_{M}}{dt}\left(  t\right)  \leq0.
\]
The argument for $f_{m}\left(  t\right)  $ is similar. \qed
\end{pf}

We may find operators that are of the form
\[
\sum_{\left\{  i,j\right\}  \in\mathcal{S}_{1}}a_{ij}\left(  t\right)  \left(
f_{j}\left(  t\right)  -f_{i}\left(  t\right)  \right)
\]
but some coefficients are negative. The argument above does not work, but it
is possible that a maximum principle still holds if the sum is positive when
$f_{i}$ is minimal and the sum is negative when $f_{i}$ is maximal even though
each term is not positive or negative respectively. This motivates our
definition of parabolic-like operators as operators which satisfy the maximum
principle for some function.

\begin{defn}
\label{parabolic-like def}An operator
\[
P:\mathcal{F}\rightarrow\mathcal{F}%
\]
of the form%
\[
\left(  Pf\right)  _{i}\left(  t\right)  =\frac{df_{i}\left(  t\right)  }%
{dt}-\sum_{\left\{  i,j\right\}  \in\mathcal{S}_{1}}a_{ij}\left(  t\right)
\left(  f_{j}\left(  t\right)  -f_{i}\left(  t\right)  \right)  ,
\]
where $a_{ij}:\left[  A,Z\right)  \rightarrow\mathbb{R}$ are self-adjoint with
respect to some metrics $\left\{  b_{i}\left(  t\right)  \right\}
_{i\in\mathcal{S}_{0}},$ is called \emph{parabolic-like }for a function
$g$\emph{ }if $Pg=0$ implies%
\begin{align*}
\frac{dg_{M}}{dt} &  \leq0\\
\frac{dg_{m}}{dt} &  \geq0
\end{align*}
for all $t\in\left[  A,Z\right)  .$
\end{defn}

Parabolic-like operators formally look the same as parabolic operators, but
some of the coefficients $a_{ij}$ may be negative. We shall show that our
discrete curvature flow equation is parabolic-like for the curvature function
in a large subset of the domain. The hope is that this will be enough to prove
pinching and convergence theorems.

We note that the maximum and minimum may be done separately, defining
\emph{upper parabolic-like} and \emph{lower parabolic-like} in the obvious
ways. In some situations we may be interested in only one of these.

\section{Combinatorial Yamabe flow}

Here we reintroduce the concepts of combinatorial Yamabe flow, based on the
work on the combinatorial Ricci flow in two dimensions by Chow-Luo
\cite{chowluo} and the combinatorial scalar curvature by Cooper-Rivin
\cite{cooperrivin}. Further details can be found in
\cite{glickensteincombinatorialyamabeflow}. Recall that if $\mathcal{S=}%
\left\{  \mathcal{S}_{0},\mathcal{S}_{1},\mathcal{S}_{2},\mathcal{S}%
_{3}\right\}  $ is a simplicial complex of dimension $3,$ where $\mathcal{S}%
_{i}$ is the $i$-dimensional skeleton, we define the metric structure as a map%
\[
r:\mathcal{S}_{0}\mathcal{\rightarrow}\left(  0,\infty\right)
\]
such that for every edge $\left\{  i,j\right\}  \in\mathcal{S}_{1}$ between
vertices $i$ and $j,$ the length of the edge is $\ell_{ij}=r_{i}+r_{j}.$ The
set of all such metrics is called the conformal class since rescaling the
$r_{i}$ will deform the structure to the metric structure with all edges the
same length. We shall use $\mathcal{T=}\left\{  \mathcal{T}_{0},\mathcal{T}%
_{1},\mathcal{T}_{2},\mathcal{T}_{3}\right\}  $ to denote the triangulation of
one tetrahedron. Recall that in order for 4 positive numbers $r_{i}%
,r_{j},r_{k},r_{\ell}$ to define a nondegenerate tetrahedron, they must
satisfy the Descartes inequality%
\[
Q_{ijk\ell}=\left(  \frac{1}{r_{i}}+\frac{1}{r_{j}}+\frac{1}{r_{k}}+\frac
{1}{r_{\ell}}\right)  ^{2}-2\left(  \frac{1}{r_{i}^{2}}+\frac{1}{r_{j}^{2}%
}+\frac{1}{r_{k}^{2}}+\frac{1}{r_{\ell}^{2}}\right)  >0.
\]
We refer to $Q_{ijk\ell}$ as the nondegeneracy quadratic. $Q_{ijk\ell}$ is
related both to the volume of the tetrahedron and the radius of the
circumscripted sphere, i.e. the sphere which is tangent to all six edges of
the tetrahedron. For a given tetrahedron, the existence of a circumscripted
sphere is equivalent to being able to define the lengths by assigning weights
$r_{i}$ to the vertices.

The curvature $K_{i}$ associated to a vertex $i$ is defined as
\[
K_{i}\doteqdot4\pi-\sum_{\left\{  i,j,k,\ell\right\}  \in\mathcal{S}_{3}%
}\alpha_{ijk\ell}%
\]
where $\alpha_{ijk\ell}$ is the solid angle vertex $i$ in the tetrahedron
$\left\{  i,j,k,\ell\right\}  .$ The solid angle is the area of the triangle
on the unit sphere cut out by the planes determined by $\left\{
i,j,k\right\}  ,\left\{  i,j,\ell\right\}  ,\left\{  i,k,\ell\right\}  $ where
$i$ is the center of the sphere. Note that $\alpha_{ijk\ell}$ is symmetric in
all permutations of the last three indices; when it is clear which tetrahedron
we are working with, we will use the simplified notation $\alpha_{i}.$ The
combinatorial Yamabe flow is defined to be
\[
\frac{dr_{i}}{dt}=-K_{i}r_{i}.
\]
Careful calculation shows that curvature satisfies the following evolution%
\[
\frac{dK_{i}}{dt}=\sum_{\left\{  i,j,k,\ell\right\}  \in\mathcal{S}_{3}%
}\left[  \Omega_{ijk\ell}\left(  K_{j}-K_{i}\right)  +\Omega_{ikj\ell}\left(
K_{k}-K_{i}\right)  +\Omega_{i\ell jk}\left(  K_{\ell}-K_{i}\right)  \right]
.
\]
This form is gotten by using the Schl\"{a}fli formula (see, for instance,
\cite{milnorschlafli}) which can be written in the following way in this case:%
\begin{equation}
r_{i}\frac{\partial\alpha_{ijk\ell}}{\partial r_{i}}+r_{j}\frac{\partial
\alpha_{ijk\ell}}{\partial r_{j}}+r_{k}\frac{\partial\alpha_{ijk\ell}%
}{\partial r_{k}}+r_{\ell}\frac{\partial\alpha_{ijk\ell}}{\partial r_{\ell}%
}=0.\label{schlafli}%
\end{equation}
We computed the partial derivatives of the angles to be
\begin{align}
\frac{\partial\alpha_{ijk\ell}}{\partial r_{i}} &  =-\frac{8r_{j}^{2}r_{k}%
^{2}r_{\ell}^{2}}{3P_{ijk}P_{ij\ell}P_{ik\ell}V_{ijk\ell}}\left[  \left(
\frac{2}{r_{i}}+\frac{1}{r_{j}}+\frac{1}{r_{k}}+\frac{1}{r_{\ell}}\right)
\right.  \label{partialairi}\\
&  \quad\quad+\frac{r_{j}}{r_{i}}\left(  \frac{1}{r_{i}}+\frac{1}{r_{k}}%
+\frac{1}{r_{\ell}}\right)  +\frac{r_{k}}{r_{i}}\left(  \frac{1}{r_{i}}%
+\frac{1}{r_{j}}+\frac{1}{r_{\ell}}\right)  \nonumber\\
&  \quad\quad+\left.  \frac{r_{\ell}}{r_{i}}\left(  \frac{1}{r_{i}}+\frac
{1}{r_{j}}+\frac{1}{r_{k}}\right)  +\left(  2r_{i}+r_{j}+r_{k}+r_{\ell
}\right)  Q_{ijk\ell}\right]  \nonumber
\end{align}
and
\begin{align}
\frac{\partial\alpha_{ijk\ell}}{\partial r_{j}} &  =\frac{4r_{i}r_{j}r_{k}%
^{2}r_{\ell}^{2}}{3P_{ijk}P_{ij\ell}V_{ijk\ell}}\left(  \frac{1}{r_{i}}\left(
\frac{1}{r_{j}}+\frac{1}{r_{k}}+\frac{1}{r_{\ell}}\right)  +\frac{1}{r_{j}%
}\left(  \frac{1}{r_{i}}+\frac{1}{r_{k}}+\frac{1}{r_{\ell}}\right)  \right.
\label{partialairj}\\
&  \quad\quad\quad\quad-\left.  \left(  \frac{1}{r_{k}}-\frac{1}{r_{\ell}%
}\right)  ^{2}\right)  ,\nonumber
\end{align}
where $P_{ijk}$ is the perimeter of the triangle $\left\{  i,j,k\right\}  $
and $V_{ijk\ell}$ is the volume of the tetrahedron $\left\{  i,j,k,\ell
\right\}  .$ The coefficients $\Omega_{ijk\ell}$ are
\[
\Omega_{ijk\ell}=\frac{\partial\alpha_{ijk\ell}}{\partial r_{j}}r_{j}%
\]
and are thus easily computed to be
\begin{align}
\Omega_{ijk\ell} &  =\frac{4r_{i}r_{j}^{2}r_{k}^{2}r_{\ell}^{2}}%
{3P_{ijk}P_{ij\ell}V_{ijk\ell}}\left(  \frac{1}{r_{i}}\left(  \frac{1}{r_{j}%
}+\frac{1}{r_{k}}+\frac{1}{r_{\ell}}\right)  +\frac{1}{r_{j}}\left(  \frac
{1}{r_{i}}+\frac{1}{r_{k}}+\frac{1}{r_{\ell}}\right)  \right.  \label{coeffs}%
\\
&  \quad\quad\quad\quad-\left.  \left(  \frac{1}{r_{k}}-\frac{1}{r_{\ell}%
}\right)  ^{2}\right)  .\nonumber
\end{align}
Geometrically, we find that%
\[
\sum_{\left\{  i,j,k,\ell\right\}  \in\mathcal{S}_{3}}\Omega_{ijk\ell}%
=\frac{\ell_{ij}^{\ast}}{r_{i}\ell_{ij}}%
\]
where $\ell_{ij}^{\ast}$ is the area of the dual faces and the sum is over all
$k$ and $\ell$, where duality comes from assigning the geometric dual to a
tetrahedron to be the center of the circumscripted sphere. The evolution of
curvature can be written compactly as
\[
\frac{dK_{i}}{dt}=\triangle K_{i}%
\]
if we define the operator $\triangle$ as
\begin{equation}
\triangle f_{i}=\frac{1}{r_{i}}\sum_{\left\{  i,j\right\}  \in\mathcal{S}_{1}%
}\frac{\ell_{ij}^{\ast}}{\ell_{ij}}\left(  f_{j}-f_{i}\right)
.\label{laplacebeltrami}%
\end{equation}
Note that this operator looks like%
\[
\triangle f_{i}=\sum_{\left\{  i,j\right\}  \in\mathcal{S}_{1}}a_{ij}\left(
f_{j}-f_{i}\right)  ,
\]
which is a Laplacian on the graph $\left(  \mathcal{S}_{0},\mathcal{S}%
_{1}\right)  $ with weights $a_{ij},$ except that it is possible for $a_{ij}$
to be negative. It can be argued that $\triangle$ is a discrete analogue of
the Laplace-Beltrami operator. We shall see that since $\Omega_{ijk\ell}$ are
not always positive, the maximum principle is not ensured. In the next section
we explore the maximum principle for this Laplacian.

\section{Parabolicity of $\frac{\partial}{\partial t}-\triangle$}

The operator $\frac{\partial}{\partial t}-\triangle,$ where the
Laplace-Beltrami operator $\triangle$ is defined by (\ref{laplacebeltrami}),
may not form a parabolic operator as defined in Definition \ref{parabolic def}%
. This is because the coefficients $\Omega_{ijk\ell}=\frac{\partial\alpha_{i}%
}{\partial r_{j}}r_{j}$ may, in fact be negative. We can see this by using our
explicit calculation for the coefficients in (\ref{coeffs}). If we choose, for
instance, $r_{1}=r_{2}=r_{3}=1$ and $r_{4}=1/5$ then we see that the
tetrahedron is not degenerate since
\[
Q_{1234}=8>0
\]
and that one coefficient is negative,%

\begin{equation}
\Omega_{1234}=-\frac{8}{75P_{123}P_{124}V_{1234}}<0.\nonumber
\end{equation}
This indicates that the maximum principle will not hold in general. However,
there are several indications that this operator is a good operator and that
it might be parabolic-like for the curvature in the sense of Definition
\ref{parabolic-like def}. First is that the matrix $\frac{\partial\alpha_{i}%
}{\partial r_{j}}$ is negative semidefinite, where the nullspace is spanned by
the vector $\left(  r_{1},r_{2},r_{3},r_{4}\right)  .$ This is shown in
Appendix A. It was originally stated in \cite{cooperrivin} but the proof there
is incorrect. It is interesting to note that the incorrect proof tried to show
that the coefficients of the Laplace-Beltrami operator are all positive, which
we have shown to be false by a difficult calculation. Clearly the case when
the coefficients are all positive is a large set.

We pursue a different understanding of why the maximum principle should hold.
The first evidence is that the Schl\"{a}fli formula tells us that
\[
r_{i}\frac{\partial\alpha_{ijk\ell}}{\partial r_{i}}+r_{j}\frac{\partial
\alpha_{jik\ell}}{\partial r_{i}}+r_{k}\frac{\partial\alpha_{kij\ell}%
}{\partial r_{i}}+r_{\ell}\frac{\partial\alpha_{\ell ijk}}{\partial r_{i}}=0
\]
while the formula (\ref{partialairi}) for $\frac{\partial\alpha_{ijk\ell}%
}{\partial r_{i}}$ calculated in \cite{glickensteincombinatorialyamabeflow}
shows that
\[
\frac{\partial\alpha_{ijk\ell}}{\partial r_{i}}<0
\]
for any tetrahedron. Together with the following lemma, we shall see that
there are large restrictions on when coefficients $\Omega_{ijk\ell}$ can be negative.

\begin{lem}
\label{negativeomega}If $\Omega_{abcd}<0$ then $r_{a},r_{b}>\min\left\{
r_{c},r_{d}\right\}  .$
\end{lem}

\begin{pf}
Suppose $\Omega_{abcd}<0,$ then%
\[
\frac{1}{r_{a}}\left(  \frac{1}{r_{b}}+\frac{1}{r_{c}}+\frac{1}{r_{d}}\right)
+\frac{1}{r_{b}}\left(  \frac{1}{r_{a}}+\frac{1}{r_{c}}+\frac{1}{r_{d}%
}\right)  -\left(  \frac{1}{r_{c}}-\frac{1}{r_{d}}\right)  ^{2}<0.
\]
eliminating the denominators and regrouping terms we get%
\[
r_{d}r_{c}\left[  \left(  r_{d}+r_{b}\right)  r_{c}+\left(  r_{b}%
+r_{c}\right)  r_{d}\right]  +r_{a}\left(  r_{d}r_{c}\left(  r_{d}%
+r_{c}\right)  -r_{b}\left(  r_{d}-r_{c}\right)  ^{2}\right)  <0
\]
so in particular we need
\[
r_{d}r_{c}\left(  r_{d}+r_{c}\right)  -r_{b}\left(  r_{d}-r_{c}\right)
^{2}<0.
\]
Solving for $r_{b}$ we get%
\[
r_{b}>\frac{r_{d}r_{c}\left(  r_{d}+r_{c}\right)  }{\left(  r_{d}%
-r_{c}\right)  ^{2}}\geq\frac{r_{d}r_{c}}{\left\vert r_{d}-r_{c}\right\vert
}\geq\min\left\{  r_{d},r_{c}\right\}
\]
since $r_{d}+r_{c}\geq\left\vert r_{d}-r_{c}\right\vert $ and $\max\left\{
r_{c},r_{d}\right\}  \geq\left\vert r_{d}-r_{c}\right\vert .$ Since the
initial expression is symmetric in $a$ and $b,$ we get $r_{a}\geq\min\left\{
r_{d},r_{c}\right\}  $ too. \qed

\end{pf}

\begin{cor}
\label{sign partial derivative of angles}If $\left\{  i,j,k,\ell\right\}
\in\mathcal{S}_{3}$ and $r_{i}=\min\left\{  r_{i},r_{j},r_{k},r_{\ell
}\right\}  $ then
\[
\Omega_{ijk\ell}\leq0\text{ and }\Omega_{jik\ell}\geq0
\]
or, equivalently,%
\[
\frac{\partial\alpha_{ijk\ell}}{\partial r_{j}}\leq0\text{ and }\frac
{\partial\alpha_{jik\ell}}{\partial r_{i}}\geq0.
\]

\end{cor}

\section{Monotonicity}

We need some way to relate the coefficients $\Omega_{ijk\ell}$ and the
curvatures. We attempt this by trying to prove that the $K$'s are monotonic as
functions of the $r$'s. This will turn out to be true for the two cases of the
double tetrahedron and the boundary of a 4-simplex, though not in general. We
consider a tetrahedron determined by $\left\{  r_{1},r_{2},r_{3}%
,r_{4}\right\}  .$ First we prove a lemma about the kinds of degeneracies that
can develop.

\begin{prop}
If $Q_{ijk\ell}\rightarrow0$ without any of the $r_{i}$ going to $0,$ then one
solid angle goes to $2\pi$ and the others go to $0.$ The solid angle
$\alpha_{i}$ which goes to $2\pi$ corresponds to $r_{i}$ being the minimum.
\end{prop}

\begin{pf}
Rewrite $Q_{ijk\ell}$ as
\begin{align}
Q_{ijk\ell}  &  =\left(  \frac{1}{r_{i}}+\frac{1}{r_{j}}+\frac{1}{r_{k}}%
+\frac{1}{r_{\ell}}\right)  ^{2}-2\left(  \frac{1}{r_{i}^{2}}+\frac{1}%
{r_{j}^{2}}+\frac{1}{r_{k}^{2}}+\frac{1}{r_{\ell}^{2}}\right) \nonumber\\
&  =\frac{1}{r_{i}}\left(  \frac{1}{r_{j}}+\frac{1}{r_{k}}+\frac{1}{r_{\ell}%
}-\frac{1}{r_{i}}\right)  +\frac{1}{r_{j}}\left(  \frac{1}{r_{i}}+\frac
{1}{r_{k}}+\frac{1}{r_{\ell}}-\frac{1}{r_{j}}\right) \label{Q expanded}\\
&  \quad+\frac{1}{r_{k}}\left(  \frac{1}{r_{i}}+\frac{1}{r_{j}}+\frac
{1}{r_{\ell}}-\frac{1}{r_{k}}\right)  +\frac{1}{r_{\ell}}\left(  \frac
{1}{r_{i}}+\frac{1}{r_{j}}+\frac{1}{r_{k}}-\frac{1}{r_{\ell}}\right)
.\nonumber
\end{align}
If $r_{i}$ is the minimum, then
\begin{align*}
\frac{1}{r_{i}}+\frac{1}{r_{k}}+\frac{1}{r_{\ell}}-\frac{1}{r_{j}}  &  >0\\
\frac{1}{r_{i}}+\frac{1}{r_{j}}+\frac{1}{r_{\ell}}-\frac{1}{r_{k}}  &  >0\\
\frac{1}{r_{i}}+\frac{1}{r_{j}}+\frac{1}{r_{k}}-\frac{1}{r_{\ell}}  &  >0.
\end{align*}
Hence if $Q_{ijk\ell}=0$ then
\[
\frac{1}{r_{j}}+\frac{1}{r_{k}}+\frac{1}{r_{\ell}}-\frac{1}{r_{i}}<0.
\]
Now we look at the partial derivative
\[
\frac{\partial}{\partial r_{i}}Q_{ijk\ell}=-\frac{2}{r_{i}^{2}}\left(
\frac{1}{r_{j}}+\frac{1}{r_{k}}+\frac{1}{r_{\ell}}-\frac{1}{r_{i}}\right)
\geq0.
\]
So if $Q_{ijk\ell}=0,$ we can always increase $r_{i}$ to make the tetrahedron nondegenerate.

Now we can categorize the degenerations. Notice that by the formula for
volume,
\[
V_{ijk\ell}=\frac{2A_{ijk}A_{ij\ell}\sin\beta_{ijk\ell}}{3\ell_{ij}},
\]
we must have that if the tetrahedron degenerates, $\sin\beta_{ijk\ell
}\rightarrow0$ for all dihedral angles. Hence $\beta_{ijk\ell}$ goes to $0$ or
to $\pi.$ Since
\begin{align*}
\frac{\partial\beta_{ijk\ell}}{\partial r_{i}} &  =\frac{2r_{i}r_{j}r_{k}%
^{2}r_{\ell}^{2}}{3P_{ijk}P_{ij\ell}V_{ijk\ell}}\left[  -\frac{1}{r_{k}^{2}%
}-\frac{1}{r_{\ell}^{2}}-2\frac{r_{j}}{r_{i}}\left(  \frac{1}{r_{i}\,r_{k}%
}+\frac{1}{r_{i}\,r_{\ell}}+\frac{1}{r_{k}\,r_{\ell}}\left(  2+\frac{r_{j}%
}{r_{i}}\right)  \right)  \right.  \\
&  \quad\quad\,\left.  +\left(  \frac{1}{r_{j}}-\frac{1}{r_{i}}\right)
\left(  \frac{2}{r_{i}}+\frac{1}{r_{k}}+\frac{1}{r_{\ell}}\right)  \right]
\end{align*}
(see \cite{glickensteincombinatorialyamabeflow}), if $r_{i}$ is the minimum,
then $\frac{\partial\beta_{ijk\ell}}{\partial r_{i}}<0.$ When the tetrahedron
first becomes degenerate, we can increase $r_{i}$ to become nondegenerate
again. But this indicates that in this case $\beta_{ijk\ell}$ would decrease,
so if $\beta_{ijk\ell}=0,$ $\beta_{ijk\ell}$ would become negative in a
nondegenerate tetrahedron. This is a contradiction, so we cannot have
$\beta_{ijk\ell}=0$ in the limit. Hence
\[
\beta_{ijk\ell}=\beta_{ikj\ell}=\beta_{i\ell jk}=\pi
\]
and $\alpha_{ijk\ell}=2\pi.$ Since $\alpha_{ijk\ell}+\alpha_{jik\ell}%
+\alpha_{kij\ell}+\alpha_{\ell ijk}\leq2\pi$ in any tetrahedron (see proof in
\cite{gaddum}), we must have
\[
\alpha_{jik\ell}=\alpha_{kij\ell}=\alpha_{\ell ijk}=0.
\]
\qed
\end{pf}

Now we can prove a monotonicity formula for angles in a given simplex as follows.

\begin{lem}
\label{rs and angles}$r_{i}\leq r_{j}$ if and only if $\alpha_{ijk\ell}%
\geq\alpha_{jik\ell}.$
\end{lem}

\begin{pf}
It is equivalent to prove that the strict inequality $r_{i}<r_{j}$ implies
$\alpha_{ijk\ell}>\alpha_{jik\ell}$ and that $r_{i}=r_{j}$ implies
$\alpha_{ijk\ell}=\alpha_{jik\ell}.$ The second statement is clear. Since we
are only looking at one tetrahedron, we can use $\alpha_{i}$ instead of
$\alpha_{ijk\ell}$ without causing confusion.

Consider the path $\sigma\left(  s\right)  =\left(  \sigma_{1}\left(
s\right)  ,\sigma_{2}\left(  s\right)  ,\sigma_{3}\left(  s\right)
,\sigma_{4}\left(  s\right)  \right)  $ defined by $\sigma\left(  s\right)
=\left(  r_{i},r_{j},\left(  1-s\right)  r_{k}+sr_{\ell},sr_{k}+\left(
1-s\right)  r_{\ell}\right)  ,$ where $r_{k}<r_{\ell}.$ We can think of
$\alpha$ as a function of four variables, where $\alpha\left(  r_{i}%
,r_{j},r_{k},r_{\ell}\right)  =\alpha_{i}.$ Let $\alpha_{3}\left(
x,y,z,w\right)  =\alpha\left(  z,x,y,w\right)  $ and $\alpha_{4}\left(
x,y,z,w\right)  =\alpha\left(  w,x,y,z\right)  .$

Consider the function%
\begin{equation}
A\doteqdot\sum_{i\in\mathcal{T}_{0}}r_{i}\alpha_{i}. \label{define A}%
\end{equation}
By the Schl\"{a}fli formula, we find that%
\[
dA=\sum_{i\in\mathcal{T}_{0}}\alpha_{i}\,dr_{i}%
\]
so
\begin{align*}
\frac{\partial A}{\partial r_{i}}  &  =\alpha_{i}\\
\frac{\partial^{2}A}{\partial r_{i}\partial r_{j}}  &  =\frac{\partial
\alpha_{j}}{\partial r_{i}}=\frac{\partial\alpha_{i}}{\partial r_{j}}.
\end{align*}
More details can be found in \cite{glickensteincombinatorialyamabeflow}. Now,
consider
\[
D\left(  s\right)  \doteqdot\frac{d}{ds}A\left(  \sigma\left(  s\right)
\right)  =\left[  \alpha_{3}\left(  \sigma\left(  s\right)  \right)
-\alpha_{4}\left(  \sigma\left(  s\right)  \right)  \right]  \left(  r_{\ell
}-r_{k}\right)  .
\]
The path is constructed so that $D\left(  1/2\right)  =0.$ Since the solid
angles are between $0$ and $2\pi$, $D\left(  s\right)  \leq2\pi\left(
r_{\ell}-r_{k}\right)  $ if the tetrahedron is nondegenerate for
$\sigma\left(  s\right)  $. Consider the derivative%
\[
\frac{d}{ds}D\left(  s\right)  =\frac{d^{2}}{ds^{2}}A\left(  \sigma\left(
s\right)  \right)  =\left[  \partial_{3}\alpha_{3}\left(  \sigma\left(
s\right)  \right)  -2\partial_{4}\alpha_{3}\left(  \sigma\left(  s\right)
\right)  +\partial_{4}\alpha_{4}\left(  \sigma\left(  s\right)  \right)
\right]  \left(  r_{\ell}-r_{k}\right)  ^{2}.
\]
If the tetrahedron is nondegenerate, then $\frac{d}{ds}D\left(  s\right)  <0$
since the Hessian of $A$ is negative definite (see Appendix A) and none of the
$\sigma_{i}$ are equal to zero. Now, as we move along the path starting at
$s=0$, either the tetrahedron degenerates for some $s_{0}\in\left(
0,1/2\right)  $ or the tetrahedron is nondegenerate up to $s=1/2.$ Suppose the
first degeneracy is at $s=s_{0}.$ Then either
\[
\sigma_{3}\left(  s_{0}\right)  =\min_{i=1,\ldots,4}\left\{  \sigma_{i}\left(
s_{0}\right)  \right\}
\]
or not. If it is the minimum, then at the degeneracy, $\alpha_{3}=2\pi$ and
the other angles are $0.$ This cannot happen, however, since then $D\left(
s_{0}\right)  =2\pi\left(  r_{\ell}-r_{k}\right)  ,$ the maximum possible, but
the derivative is negative, since $\frac{d}{ds}D\left(  s_{0}\right)  \leq0$
implies that $D\left(  s\right)  $ must have been larger than its maximum for
some $s<s_{0},$ a contradiction. If $\sigma_{3}\left(  s\right)  $ is not the
minimum, then at a degenerate point, $\alpha_{3}=\alpha_{4}=0.$ So there
exists a first point $s_{1}\in(0,1/2]$ where $\alpha_{k}=\alpha_{\ell}$ such
that the tetrahedron is nondegenerate for $s\in\lbrack0,s_{1}).$ Hence
$D\left(  s_{1}\right)  =0.$ Since the tetrahedron is nondegenerate on
$[0,s_{1}),$ we have $\frac{d}{ds}D<0$ for $s\in\lbrack0,s_{1}).$ Together
with $D\left(  s_{1}\right)  =0$ we have $D\left(  s\right)  >0$ for
$s\in\lbrack0,s_{1}).$ In particular, $D\left(  0\right)  >0,$ i.e.
$\alpha_{k}>\alpha_{\ell}.$

The reverse inequality must be true as well since the argument is symmetric.
\qed

\end{pf}

The lemma has the following interesting geometric consequence for a conformal
tetrahedron which is an analogue of the fact that in any triangle the longest
side is opposite the largest angle. The author does not know if this statement
is true for a general tetrahedron.

\begin{cor}
For a conformal tetrahedron with metric structure $\left\{  r_{1},r_{2}%
,r_{3},r_{4}\right\}  $, the side with the largest area is opposite the
largest solid angle, and the side with the smallest angle is opposite the
smallest solid angle.
\end{cor}

\begin{pf}
This follows from the fact that the angle $\alpha_{ijk\ell}$ is opposite the
side with area $\sqrt{r_{j}r_{k}r_{\ell}\left(  r_{j}+r_{k}+r_{\ell}\right)
}.$ \qed

\end{pf}

The lemma can be used to show monotonicity of curvature for small
triangulations as follows.

\begin{cor}
\label{monotonedoubletetra}In the case of the double tetrahedron, $r_{i}\leq
r_{j}$ if and only if $K_{i}\leq K_{j}.$
\end{cor}

\begin{pf}
This follows from the fact that $K_{i}=4\pi-2\alpha_{ijk\ell}.$ \qed

\end{pf}

\begin{cor}
\label{monotone4simplex}In the case of the boundary of a $4$-simplex,
$r_{i}\leq r_{j}$ if and only if $K_{i}\leq K_{j}.$
\end{cor}

\begin{pf}
Once again it is sufficient to show that $r_{i}<r_{j}$ implies $K_{i}<K_{j}$
and $r_{i}=r_{j}$ implies $K_{i}=K_{j}.$ The latter is trivial. Let's number
the vertices $\left\{  1,\ldots,5\right\}  .$ Suppose $r_{1}<r_{2}.$ Then
\begin{align*}
K_{1}  &  =4\pi-\alpha_{1234}-\alpha_{1235}-\alpha_{1245}-\alpha_{1345}\\
K_{2}  &  =4\pi-\alpha_{2134}-\alpha_{2135}-\alpha_{2145}-\alpha_{2345}.
\end{align*}
We know by Lemma \ref{rs and angles} that $\alpha_{2134}<\alpha_{1234},$
$\alpha_{2135}<\alpha_{1235},$ and $\alpha_{2145}<\alpha_{1245}.$ We thus need
only show $\alpha_{2345}<\alpha_{1345}$. Consider the path
\[
\sigma\left(  s\right)  =\left(  \frac{r_{1}r_{2}}{\left(  1-s\right)
r_{2}+sr_{1}},r_{3},r_{4},r_{5}\right)  .
\]
We notice that that the nondegeneracy quadratic $Q\left(  \sigma\left(
s\right)  \right)  $ is a polynomial in $s$ with highest term
\[
-\left(  \frac{1}{r_{1}}-\frac{1}{r_{2}}\right)  ^{2}s^{2}.
\]
Since the quadratic $Q\left(  \sigma\left(  s\right)  \right)  $ is concave
the minimum for $s\in\left[  0,1\right]  $ must occur at $s=0$ or $s=1.$ But
since $Q\left(  \sigma\left(  0\right)  \right)  =Q_{1345}>0$ and $Q\left(
\sigma\left(  1\right)  \right)  =Q_{2345}>0,$ $Q\left(  \sigma\left(
s\right)  \right)  >0$ for all $s\in\left[  0,1\right]  .$ Furthermore
\[
\frac{d}{ds}\alpha\left(  \sigma\left(  s\right)  \right)  =\partial_{1}%
\alpha\left(  \sigma\left(  s\right)  \right)  \frac{r_{1}r_{2}\left(
r_{2}-r_{1}\right)  }{\left(  \left(  1-s\right)  r_{2}+r_{1}s\right)  ^{2}}%
\]
which is negative for all $s\in\left[  0,1\right]  $ by formula
(\ref{partialairi}) since $r_{2}>r_{1}.$ Hence $\alpha\left(  \sigma\left(
0\right)  \right)  >\alpha\left(  \sigma\left(  1\right)  \right)  ,$ or
$\alpha_{1245}>\alpha_{2145}.$ \qed

\end{pf}

Proving a similar statement about larger complexes would be much harder, since
we cannot pair up angles which are in the same or bordering tetrahedra as we
do here. We shall call this condition the \emph{monotonicity condition} for a
tetrahedron $\left\{  i,j,k,\ell\right\}  $:%
\begin{equation}
r_{i}\leq r_{j}\text{ if and only if }K_{i}\leq K_{j}. \tag{MC}%
\label{monotone}%
\end{equation}

The condition is true for an open set of triangulations. Unfortunately,
\ref{monotone} is not necessarily preserved under the flow. For instance, if
$r_{i}=r_{j}$ but $K_{i}<K_{j}$ then $\frac{dr_{i}}{dt}=-K_{i}r_{i}%
>-K_{j}r_{j}=\frac{dr_{j}}{dt}.$ The monotonicity will counteract some of the
potential degenerations of the flow, and thus allows proof of the maximum
principle and long term estimates.

\section{Proof of the maximum principle}

Suppose we have a complex such that each tetrahedron $\left\{  i,j,k,\ell
\right\}  $ satisfies the monotonicity condition \ref{monotone}. Assume that
$r_{i}\leq r_{j},r_{k}\leq r_{\ell}$. We shall first look at the minimum. By
Corollary \ref{sign partial derivative of angles}, $\Omega_{ijk\ell}%
,\Omega_{ikj\ell},\Omega_{i\ell jk}$ are all nonnegative, so since $K_{i}$ is
the minimum curvature among $\left\{  i,j,k,\ell\right\}  $,
\[
\Omega_{ijk\ell}\left(  K_{j}-K_{i}\right)  +\Omega_{ikj\ell}\left(
K_{k}-K_{i}\right)  +\Omega_{i\ell jk}\left(  K_{\ell}-K_{i}\right)  \geq0.
\]
Now for the entire triangulation, we get that if $K_{m}=\min_{i\in
\mathcal{S}_{0}}\left\{  K_{i}\right\}  ,$ then
\[
\frac{d}{dt}K_{m}=\sum_{\left\{  m,j,k,\ell\right\}  \in\mathcal{S}_{3}%
}\left[  \Omega_{mjk\ell}\left(  K_{j}-K_{m}\right)  +\Omega_{mkj\ell}\left(
K_{k}-K_{m}\right)  +\Omega_{i\ell jk}\left(  K_{\ell}-K_{m}\right)  \right]
\]
must be a sum of nonnegative numbers since in any tetrahedron containing $m,$
$r_{i}=r_{m}$ is the smallest weight.

We now look at the maximum. We want to show that
\[
\Omega_{\ell ijk}\left(  K_{i}-K_{\ell}\right)  +\Omega_{\ell jik}\left(
K_{j}-K_{\ell}\right)  +\Omega_{\ell kij}\left(  K_{k}-K_{\ell}\right)
\leq0.
\]
This is certainly true if $\Omega_{\ell ijk},\Omega_{\ell jik},\Omega_{\ell
kij}$ are all nonnegative since $K_{\ell}$ is the largest curvature. Again
using Lemma \ref{negativeomega} we see that if $\Omega_{abcd}<0$ then $b\neq
i,$ so $\Omega_{\ell ijk}\geq0.$ We are then left with the case that both
$\Omega_{\ell jik},\Omega_{\ell kij}$ are negative or only one is negative.
First consider the case when both are negative. In this case it is sufficient
to show that
\[
\Omega_{\ell ijk}+\Omega_{\ell jik}+\Omega_{\ell kij}\geq0
\]
since in this case we have
\begin{align*}
\Omega_{\ell ijk}  &  \geq-\Omega_{\ell jik}-\Omega_{\ell kij}\\
\Omega_{\ell ijk}\left(  K_{i}-K_{\ell}\right)   &  \leq-\left(  \Omega_{\ell
jik}+\Omega_{\ell kij}\right)  \min\left\{  \left(  K_{j}-K_{\ell}\right)
,\left(  K_{k}-K_{\ell}\right)  \right\} \\
&  \leq-\Omega_{\ell jik}\left(  K_{j}-K_{\ell}\right)  -\Omega_{\ell
kij}\left(  K_{k}-K_{\ell}\right)
\end{align*}
since
\begin{align*}
\left(  K_{i}-K_{\ell}\right)   &  \leq\min\left\{  \left(  K_{j}-K_{\ell
}\right)  ,\left(  K_{k}-K_{\ell}\right)  \right\} \\
&  \leq\max\left\{  \left(  K_{j}-K_{\ell}\right)  ,\left(  K_{k}-K_{\ell
}\right)  \right\}
\end{align*}
which is nonpositive and $\Omega_{\ell ijk},-\Omega_{\ell jik},-\Omega_{\ell
kij}\geq0.$ So
\[
\Omega_{\ell ijk}\left(  K_{i}-K_{\ell}\right)  +\Omega_{\ell jik}\left(
K_{j}-K_{\ell}\right)  +\Omega_{\ell kij}\left(  K_{k}-K_{\ell}\right)  \leq0
\]
and we are done. The inequality $\Omega_{\ell ijk}+\Omega_{\ell jik}%
+\Omega_{\ell kij}\geq0$ follows from the Schl\"{a}fli formula, since%
\begin{align*}
\Omega_{\ell ijk}+\Omega_{\ell jik}+\Omega_{\ell kij}  &  =\frac
{\partial\alpha_{\ell ijk}}{\partial r_{i}}r_{i}+\frac{\partial\alpha_{\ell
ijk}}{\partial r_{j}}r_{j}+\frac{\partial\alpha_{\ell ijk}}{\partial r_{k}%
}r_{k}\\
&  =-\frac{\partial\alpha_{\ell ijk}}{\partial r_{\ell}}r_{\ell}.
\end{align*}
This is nonnegative by formula (\ref{partialairi}) for $\frac{\partial
\alpha_{\ell ijk}}{\partial r_{\ell}}.$

Now suppose that only one is negative, say $\Omega_{\ell kij}<0,$ then
similarly it is enough to show
\[
\Omega_{\ell ijk}+\Omega_{\ell kij}\geq0
\]
since then
\begin{align*}
\Omega_{\ell ijk}  &  \geq-\Omega_{\ell kij}\\
\Omega_{\ell ijk}\left(  K_{i}-K_{\ell}\right)   &  \leq-\Omega_{\ell
kij}\left(  K_{k}-K_{\ell}\right)
\end{align*}
and $\Omega_{\ell jik}\left(  K_{j}-K_{\ell}\right)  \leq0.$ We argue from our
explicit calculations. The sum $\Omega_{\ell ijk}+\Omega_{\ell kij}$ is equal to%

\begin{align*}
&  \frac{4r_{\ell}r_{j}^{2}r_{k}^{2}r_{i}^{2}}{3P_{\ell ik}P_{\ell
ij}V_{ijk\ell}}\left(  \frac{1}{r_{\ell}}\left(  \frac{1}{r_{j}}+\frac
{1}{r_{k}}+\frac{1}{r_{i}}\right)  +\frac{1}{r_{i}}\left(  \frac{1}{r_{\ell}%
}+\frac{1}{r_{j}}+\allowbreak\frac{1}{r_{k}}\right)  -\left(  \frac{1}{r_{j}%
}-\frac{1}{r_{k}}\right)  ^{2}\right) \\
&  \quad+\frac{4r_{\ell}r_{j}^{2}r_{k}^{2}r_{i}^{2}}{3P_{\ell jk}P_{\ell
ik}V_{ijk\ell}}\left(  \frac{1}{r_{\ell}}\left(  \frac{1}{r_{j}}+\frac
{1}{r_{k}}+\frac{1}{r_{i}}\right)  +\frac{1}{r_{k}}\left(  \frac{1}{r_{\ell}%
}+\frac{1}{r_{j}}+\allowbreak\frac{1}{r_{i}}\right)  -\left(  \frac{1}{r_{j}%
}-\frac{1}{r_{i}}\right)  ^{2}\right)
\end{align*}
which, when simplifying and multiplying by a positive factor, has the same
sign as
\begin{align*}
&  \left(  P_{\ell jk}+2P_{\ell ij}\right)  \frac{1}{r_{i}r_{j}}+\left(
P_{\ell jk}+\allowbreak P_{\ell ij}\right)  \frac{1}{r_{i}r_{k}}+\left(
P_{\ell ij}+2P_{\ell jk}\right)  \frac{1}{r_{\ell}r_{i}}\\
&  \quad+\allowbreak\left(  2P_{\ell jk}+P_{\ell ij}\right)  \frac{1}%
{r_{j}r_{k}}+\left(  P_{\ell ij}+P_{\ell jk}\right)  \frac{1}{r_{\ell}r_{j}%
}+\left(  2P_{\ell ij}+P_{\ell jk}\right)  \frac{1}{r_{\ell}r_{k}}\\
&  \quad-P_{\ell ij}\frac{1}{r_{i}^{2}}-\left(  P_{\ell jk}+P_{\ell
ij}\right)  \frac{1}{r_{j}^{2}}-P_{\ell jk}\frac{1}{r_{k}^{2}}.
\end{align*}
This is greater than or equal to
\[
P_{\ell ij}Q_{ijk\ell}+P_{\ell jk}\left(  \frac{1}{r_{i}r_{j}}-\frac{1}%
{r_{j}^{2}}\right)  +\left(  P_{\ell jk}-P_{\ell ij}\right)  \left(  \frac
{1}{r_{i}r_{k}}-\frac{1}{r_{k}^{2}}\right)
\]
which is nonnegative since $r_{i}\leq r_{j}$ and $\left\{  i,j,k,\ell\right\}
$ is nondegenerate.

Thus we have proven the following:

\begin{thm}
\label{maxprinc}On a complex with a metric structure which satisfies the
monotonicity condition \ref{monotone}, if $K_{m}$ is the minimum curvature,
$K_{M}$ is the maximum curvature, and we satisfy the combinatorial Yamabe
flow
\[
\frac{d}{dt}r_{i}=-K_{i}r_{i}%
\]
then
\[
\frac{d}{dt}K_{m}\geq0
\]
and
\[
\frac{d}{dt}K_{M}\leq0,
\]
i.e. the combinatorial Yamabe flow is parabolic-like for $K.$
\end{thm}

\begin{cor}
On the double tetrahedron and the boundary of a 4-simplex the combinatorial
Yamabe flow is parabolic-like for the curvature function $K$, i.e.%
\begin{align*}
\frac{d}{dt}K_{M}  &  \leq0\\
\frac{d}{dt}K_{m}  &  \geq0.
\end{align*}

\end{cor}

\begin{pf}
This follows from Corollaries \ref{monotonedoubletetra} and
\ref{monotone4simplex} which say that the monotonicity condition is satisfied
in these cases. \qed

\end{pf}

The maximum principle has been used by Hamilton and others to prove many
pinching results for geometric evolution equations (e.g.
\cite{hamilton:threemanifolds}, \cite{hamilton:singularities}). The most basic
use is to show preservation of positive or negative curvature, which is an
easy corollary.

\begin{cor}
On a complex with a metric structure which satisfies the monotonicity
condition \ref{monotone} for $t\in\lbrack0,T)$, then nonnegative curvature is
preserved, i.e. if $K_{i}\geq0$ for all $i$ for $t=0,$ then $K_{i}\geq0$ for
all $i$ for all $t\in\lbrack0,T).$ Similarly, nonpositive curvature is preserved.
\end{cor}

\begin{pf}
If $K_{i}\geq0$ for all $i,$ then in particular the minimum is nonnegative and
increasing. \qed

\end{pf}

\section{Long term existence}

The monotonicity condition (\ref{monotone}) also gives us a way to show long
term existence. The maximum principle will allow us to bound the growth or
decay of the lengths $r_{i}.$

In order to show long term existence, we must show that $Q_{ijk\ell}>0$ for
every $\left\{  i,j,k,\ell\right\}  \in\mathcal{S}_{3}$ and that $r_{i}$ does
not go to zero or infinity in finite time for any $i\in\mathcal{S}_{0}.$ Since
we are only working with one tetrahedron, we can use $Q=Q_{ijk\ell}$ without
fear of confusion. We calculate%
\begin{equation}
\frac{\partial Q}{\partial r_{i}}=-\frac{2}{r_{i}^{2}}\left(  \frac{1}{r_{j}%
}+\frac{1}{r_{k}}+\frac{1}{r_{\ell}}-\frac{1}{r_{i}}\right)  \label{partial q}%
\end{equation}
and hence using formula (\ref{Q expanded}) we see%
\[
2Q=-\left(  \frac{\partial Q}{\partial r_{i}}r_{i}+\frac{\partial Q}{\partial
r_{j}}r_{j}+\frac{\partial Q}{\partial r_{k}}r_{k}+\frac{\partial Q}{\partial
r_{\ell}}r_{\ell}\right)  .
\]
If $Q=0,$ then we have
\begin{equation}
\frac{\partial Q}{\partial r_{i}}r_{i}+\frac{\partial Q}{\partial r_{j}}%
r_{j}+\frac{\partial Q}{\partial r_{k}}r_{k}+\frac{\partial Q}{\partial
r_{\ell}}r_{\ell}=0. \label{Q=0 simplification}%
\end{equation}
Write the evolution of $Q$ as%
\begin{align*}
\frac{dQ}{dt}  &  =-\frac{\partial Q}{\partial r_{i}}K_{i}r_{i}-\frac{\partial
Q}{\partial r_{j}}K_{j}r_{j}-\frac{\partial Q}{\partial r_{k}}K_{k}r_{k}%
-\frac{\partial Q}{\partial r_{\ell}}K_{\ell}r_{\ell}\\
&  =-\frac{\partial Q}{\partial r_{j}}r_{j}\left(  K_{j}-K_{i}\right)
-\frac{\partial Q}{\partial r_{k}}r_{k}\left(  K_{k}-K_{i}\right)
-\frac{\partial Q}{\partial r_{\ell}}r_{\ell}\left(  K_{\ell}-K_{i}\right)
\end{align*}
using (\ref{Q=0 simplification}). Notice that if $r_{i}$ is the minimum, then
$\frac{\partial Q}{\partial r_{j}}\leq0$ for $j\neq i$ by (\ref{partial q}).
If the monotonicity condition \ref{monotone} is satisfied, then we must have%
\[
\frac{dQ}{dt}\geq0
\]
if $Q=0$ since $K_{j}\geq K_{i}$ and hence the tetrahedra do not degenerate.

Now we can show that the maximum principle gives bounds on growth and decay of
the $r_{i}$.

\begin{prop}
If $\frac{dK_{M}}{dt}\leq0$ and $\frac{dK_{m}}{dt}\geq0$ and then there is a
constant $C$ such that
\[
r_{i}\left(  0\right)  e^{-Ct}\leq r_{i}\left(  t\right)  \leq r_{i}\left(
0\right)  e^{Ct}.
\]

\end{prop}

\begin{pf}
Let $C\left(  t\right)  =\max\left\{  K_{M}\left(  t\right)  ,-K_{m}\left(
t\right)  \right\}  \geq0.$ Then
\[
-Cr_{i}\leq-K_{i}r_{i}\leq Cr_{i}%
\]
for each $t.$ Since $\frac{dK_{M}}{dt}\leq0$ and $\frac{dK_{m}}{dt}\geq0,$ we
must have $C\left(  t\right)  \leq C\left(  0\right)  .$ So if we look at the
evolution%
\[
-C\left(  0\right)  r_{i}\leq\frac{dr_{i}}{dt}\leq C\left(  0\right)  r_{i}%
\]
we get
\[
r_{i}\left(  0\right)  e^{-C\left(  0\right)  t}\leq r_{i}\left(  t\right)
\leq r_{i}\left(  0\right)  e^{C\left(  0\right)  t}.
\]
\qed

\end{pf}

Thus the solution exists for all time. Convergence to constant curvature now
follows from \cite{glickensteincombinatorialyamabeflow}.

\section{Further remarks}

In this paper we have seen two large sets of possible metric structures within
which the maximum principle holds: the set where the coefficients
$\Omega_{ijk\ell}$ are positive and the set where the monotonicity condition
\ref{monotone} is satisfied. Unfortunately, neither of these conditions is
obviously preserved by the flow. It would be highly desirable to find a set
which is preserved by the flow within which the curvature satisfies the
maximum principle.

Numerical data suggests that the maximum principle holds in much greater
generality, even for large simplicial complexes that do not satisfy
monotonicity. Numerical simulation of the flow requires a true simplicial
complex; a CW decomposition will not work because there are not enough
vertices to allow the different tetrahedra in the complex to evolve
independently. Thus the current numerical work has been limited to certain
small triangulations (fewer than 15 vertices) of the 3-sphere, the direct
product of the 2-sphere with the circle, the twisted product of the 2-sphere
with the circle, and the 3-torus. Some of the small triangulations are due to
F.H. Lutz (see \cite{lutzthesis} and \cite{lutzmanifoldpage}). The condition
of monotonicity is not particularly well understood for large complexes
either, though it is known not to hold in general even for triangulations of
$S^{3}.$

The maximum principle is closely connected to the fact that the operator
$\triangle$ is negative semi-definite in the smooth case, but it is not clear
that the definition of maximum principle for graph Laplacians which we use
here is the right maximum principle to correspond to the definiteness. This
may also be related to the fact that the principle eigenfunction in the smooth
case is positive, while in the discrete case here the principle eigenvector
usually will not have all positive (or all negative) entries. Perhaps a kind
of `refined maximum principle' is needed, as in
\cite{berestyckinirenbergvaradhan} and \cite{padilla}. Also integral type
maximum principles have been successful in studying discrete Laplacians as in
\cite{coulhongrigoryanzucca}.

\section*{Acknowledgements}

The author would like to thank Ben Chow for introducing the combinatorial
Yamabe flow and for all his help. The author would also like to thank Feng Luo
for useful conversations and Lennie Friedlander for references on the relation
of the maximum principle and eigenvalues of the Laplacian.

\section*{Appendix A}

In this appendix we prove that the four by four symmetric matrix $\left(
\frac{\partial\alpha_{i}}{\partial r_{j}}\right)  $ has three negative
eigenvalues and one zero eigenvalue. The proof in \cite{cooperrivin} is
incorrect; Igor Rivin has given a new proof in \cite{rivincorrection}. We do
not understand part of Rivin's proof, but we can complete Rivin's arguments by
a calculation given below. The complete argument is given here.

The zero eigenvalue comes from the vector $\left(  r_{1},r_{2},r_{3}%
,r_{4}\right)  $ by the Schl\"{a}fli formula (\ref{schlafli}). Recall that the
$\left(  n-1\right)  $-dimensional minor $M\left(  i,j\right)  $ of a matrix
$M$ is the matrix with the $i$th row and the $j$th column removed. Let $A$ be
the matrix of partial derivatives, so $A_{ij}=\frac{\partial\alpha_{ijk\ell}%
}{\partial r_{j}}.$ We take as the domain of $A$ the set of $\left(
r_{1},r_{2},r_{3},r_{4}\right)  $ such that the associated tetrahedron is
nondegenerate, i.e. $r_{i}>0$ for $i=1,\ldots,4$ and $Q_{ijk\ell}>0.$

\begin{prop}
\label{compute determinant}The minors $A\left(  i,j\right)  $ have determinant%
\[
\det A\left(  i,j\right)  =\left(  -1\right)  ^{i+j+1}\frac{288V_{ijk\ell}%
}{r_{i}r_{j}P_{ijk}P_{ij\ell}P_{ik\ell}P_{jk\ell}}%
\]
and hence is nonzero if the tetrahedron $\left\{  i,j,k,\ell\right\}  $ is nondegenerate.
\end{prop}

\begin{pf}
We can do a rather lengthy calculation using (\ref{partialairi}) and
(\ref{partialairj}). Note that we need only compute the minors $A\left(
i,j\right)  $ where $i\neq j$ since $\det A=0.$ The minors $A\left(
i,i\right)  $ on the diagonal are slightly more difficult to calculate because
there are 3 entries of the more complicated form $\frac{\partial\alpha_{i}%
}{\partial r_{i}}$ instead of only 2 for the off-diagonal $A\left(
i,j\right)  $ where $i\neq j.$ \qed

\end{pf}

\begin{cor}
\label{negative semidef}The matrix $A$ is negative semidefinite, rank 3, and
the nullspace is the span of the vector $\left(  r_{1},r_{2},r_{3}%
,r_{4}\right)  .$
\end{cor}

\begin{pf}
The rank follows from Proposition \ref{compute determinant}. The nullspace
condition follows from the Schl\"{a}fli formula (\ref{schlafli}). Since the
domain is connected and the rank is always 3, the eigenvalues must always have
the same sign. We need only compute the matrix at one point, say $r_{i}=1$ for
$i=1,\ldots,4.$ The matrix $A$ at this point is easily computed to be%
\[
\frac{1}{3\sqrt{2}}\left[
\begin{array}
[c]{cccc}%
-3 & 1 & 1 & 1\\
1 & -3 & 1 & 1\\
1 & 1 & -3 & 1\\
1 & 1 & 1 & -3
\end{array}
\right]  ,
\]
which has eigenvalues $0,-\frac{2}{3}\sqrt{2},-\frac{2}{3}\sqrt{2},-\frac
{2}{3}\sqrt{2}.$ \qed

\end{pf}

\end{document}